\input amstex
\input amsppt.sty
\magnification=\magstep1
\vsize=22.2truecm
\baselineskip=16truept
%\NoBlackBoxes
%\NoRunningHeads
\nologo
\pageno=1
\topmatter
\def\Z{\Bbb Z}

\def\bg{\bigg}
\def\({\bg(}
\def\[{\bg[}
\def\){\bg)}
\def\]{\bg]}

\def\bi{\binom}

\def\bi{\binom}

\def\Proof{\noindent{\it Proof}}

\def\Ack{\medskip\noindent {\bf Acknowledgment}}
\def\pmod #1{\ (\roman{mod}\ #1)}
\def\mo #1{\ \roman{mod}\ #1}
\def\qbinom #1#2#3{\left[\matrix#1\\#2\endmatrix\right]_{#3}}

\def\cf{\Cal{F}}
\def\cg{\widehat{\Cal{F}}}
\def\cp{\Cal{P}}
\def\cq{\widehat{\Cal{P}}}
\def\rT{{\roman{T}}}
\def\jacob #1#2{\left(\frac{#1}{#2}\right)}
\def\floor #1{\left\lfloor{#1}\right\rfloor}
\title
Arithmetic properties of $q$-Fibonacci numbers and $q$-Pell numbers
\endtitle
\abstract
We investigate some arithmetic properties of the $q$-Fibonacci numbers and the $q$-Pell numbers.
\endabstract
\author
Hao Pan
\endauthor
\address
Department of Mathematics, Nanjing University,
Nanjing 210093, People's Republic of China
\endaddress
\email{haopan79\@yahoo.com.cn}\endemail
\endtopmatter
\document
\TagsOnRight
\heading
1. Introduction
\endheading
The Fibonacci numbers $F_n$ are given by
$$
F_0=0,\ F_1=1\text{ and }F_n=2F_{n-1}+F_{n-2}\text{ for }n\geq 2.
$$
For any odd prime $p$, it is well-known (cf. [7]) that
$$
F_p\equiv\jacob{5}{p}\pmod{p},\tag 1.1
$$
$$
F_{p+1}\equiv\frac{1}{2}(1+\jacob{5}{p})\pmod{p}\tag 1.2
$$
and
$$
F_{p-1}\equiv\frac{1}{2}(1-\jacob{5}{p})\pmod{p},\tag 1.3
$$
where $\jacob{\cdot}{p}$ denotes the Legendre symbol. Indeed, we have
$$
\aligned
F_p=\frac{(1+\sqrt{5})^p-(1-\sqrt{5})^p}{2^p\sqrt{5}}
\equiv\frac{(1+(\sqrt{5})^p)-(1-(\sqrt{5})^p)}{2\sqrt{5}}
=5^{(p-1)/2}\pmod{p}.
\endaligned
$$
For more results on the congruences involving the Fibonacci numbers, the readers may refer to [10], [11] and [12].

On the other hand, a sequence of polynomials $\cf_n(q)$ was firstly introduced by Schur (cf. [9]):
$$
\cf_n(q)=\cases 0\qquad&\text{if }n=0,\\1\qquad&\text{if }n=1,\\\cf_{n-1}(q)+q^{n-2}\cf_{n-2}(q)\qquad&\text{if }n\geq2.\endcases
$$
Also Schur considered another sequence $\cg_n(q)$, which is given by
$$
\cg_n(q)=\cases 0\qquad&\text{if }n=0,\\1\qquad&\text{if }n=1,\\\cg_{n-1}(q)+q^{n-1}\cg_{n-2}(q)\qquad&\text{if }n\geq2.\endcases
$$
Obviously both $\cf_n(q)$ and $\cg_n(q)$ are the $q$-analogues of the Fibonacci numbers. The sequences $\cf_n(q)$ and $\cg_n(q)$
have been investigated in several papers (e.g., see [2], [3], [4], [5] and [6]).
However, seemingly there are no simple expressions for $\cf_n(q)$ and $\cg_n(q)$.

Now we can give the $q$-analogues of (1.1), (1.2) and
(1.3). Suppose that $p\not=5$ is an odd prime. Let $\alpha_p$ be the integer such that $1\leq \alpha_p\leq 4$ and
$\alpha_pp\equiv 1\mo{5}$. As usual we set
$$
[n]_q=\frac{1-q^n}{1-q}=1+q+\cdots+q^{n-1}
$$
for any non-negative integer $n$.
\proclaim{Theorem 1.1}
$$
\cf_{p+1}\equiv \frac{1}{2}(1+\jacob{5}{p})\pmod{[p]_q}\tag 1.4
$$
and
$$
\cf_{p}\equiv\jacob{5}{p}q^{((5-\alpha_p)p+1)/5}\pmod{[p]_q}.\tag 1.5
$$
\endproclaim
\proclaim{Theorem 1.2}
$$
\cg_{p-1}\equiv \frac{1}{2}(1-\jacob{5}{p})\pmod{[p]_q}\tag 1.6
$$
and
$$
\cg_{p}\equiv\jacob{5}{p}q^{(\alpha_pp-1)/5}\pmod{[p]_q}.\tag 1.7
$$
\endproclaim

The Pell numbers $P_n$ are given by
$$
P_0=0,\ P_1=1\text{ and }P_n=2P_{n-1}+P_{n-2}\text{ for }n\geq 2.
$$
It is easy to check that
$$
P_n=\frac{(1+\sqrt{2})^n-(1-\sqrt{2})^n}{2\sqrt{2}}.
$$
Hence for odd prime $p$, we have
$$
P_p=\frac{(1+\sqrt{2})^p-(1-\sqrt{2})^p}{2\sqrt{2}}\equiv\frac{2(\sqrt{2})^p}{2\sqrt{2}}=2^{(p-1)/2}\equiv\jacob{2}{p}\pmod{p}.\tag 1.8
$$
Define the $q$-Pell numbers $\cp_n(q)$ and $\cq_n(q)$ by
$$
\cp_n(q)=\cases 0\qquad&\text{if }n=0,\\1\qquad&\text{if }n=1,\\(1+q^{n-1})\cp_{n-1}(q)+q^{n-2}\cp_{n-2}(q)\qquad&\text{if }n\geq2,\endcases
$$
and
$$
\cq_n(q)=\cases 0\qquad&\text{if }n=0,\\1\qquad&\text{if }n=1,\\(1+q^{n-1})\cq_{n-1}(q)+q^{n-1}\cq_{n-2}(q)\qquad&\text{if }n\geq2.\endcases
$$
\proclaim{Theorem 1.3} Let $p$ be an odd prime. Then
$$
q^{(p^2-1)/8}\cp_p(q)\equiv\jacob{2}{p}\pmod{[p]_q}\tag 1.9
$$
and
$$
\cq_p(q)\equiv\jacob{2}{p}q^{(p^2-1)/8}\pmod{[p]_q}.\tag 1.10
$$
Furthermore, we have
$$
\cp_{p+1}(q)-\cp_{p}(q)\equiv\cq_{p+1}(q)-\cq_{p}(q)\equiv 1\pmod{[p]_q}.\tag 1.11
$$
\endproclaim
Since $\jacob{2}{p}=(-1)^{(p^2-1)/8}$, (1.9) and (1.10) can be respectively rewritten as
$$
(-q)^{(p^2-1)/8}\cp_p(q)\equiv1\pmod{[p]_q}
$$
and
$$
\cq_p(q)\equiv(-q)^{(p^2-1)/8}\pmod{[p]_q}.
$$

The proofs of Theorems 1.1, 1.2 and 1.3 will be given in Sections 2 and 3.

\heading
2. Proofs of Theorems 1.1 and 1.2
\endheading
For any $n,m\in\Z$, the $q$-binomial coefficient $\qbinom{n}{m}{q}$ is given by
$$
\qbinom{n}{m}{q}=\frac{(1-q^{n})(1-q^{n-1})\cdots(1-q^{n-m+1})}{(1-q^m)(1-q^{m-1})\cdots(1-q)}
$$
when $m\geq 0$, and let $\qbinom{n}{m}{q}=0$ if $m<0$. Obviously $\qbinom{n}{m}{q}$ is a polynomial in
the variable $q$ with integral coefficients since the $q$-binomial coefficients satisfy the recurrence relation
$$
\qbinom{n+1}{m}{q}=q^m\qbinom{n}{m}{q}+\qbinom{n}{m-1}{q}.
$$

Let $\floor{x}$ denotes the greatest integer not exceeding $x$.
Then for any non-negative integer $n$, we have
$$
\align
\cf_{n+1}(q)=&\sum_{0\leq 2j\leq n}q^{j^2}\qbinom{n-j}{j}{q}\\
=&\sum_{j=-\infty}^{\infty}(-1)^jq^{j(5j+1)/2}\qbinom{n}{\lfloor(n-5j)/2\rfloor}{q}\tag 2.1
\endalign
$$
and
$$
\align
\cg_{n+1}(q)=&\sum_{0\leq 2j\leq n}q^{j^2+j}\qbinom{n-j}{j}{q}\\
=&\sum_{j=-\infty}^{\infty}(-1)^jq^{j(5j-3)/2}\qbinom{n+1}{\lfloor(n+1-5j)/2\rfloor+1}{q}.\tag
2.2
\endalign
$$
We mention that (2.1) and (2.2) can be considered as the finite
forms of the first and the second of the Rogers-Ramanujan
identities respectively (the full proofs of (2.1) and (2.2) can
been found in [1]).

In this section we assume that $p\not=5$ is an odd prime.
\proclaim{Lemma 2.1} Let
$$
L(j)=\frac{j(5j+1)}{2}-\bi{\floor{(p-1-5j)/2}+1}{2},
$$
and let
$$
\widehat{L}(j)=\frac{j(5j-3)}{2}-\bi{\floor{p-1-5j)/2}+2}{2}.
$$
Then
$$
L(2j)-L(2j-1)=\widehat{L}(2j)-\widehat{L}(2j-1)=p.
$$
\endproclaim
\proclaim{Lemma 2.2} Let
$$
S_p=\{j\in\Z:\,0\leq\floor{(p-1-5j)/2}\leq p-1\}
$$
and
$$
\widehat{S}_p=\{j\in\Z:\,0\leq\floor{(p-1-5j)/2}+1\leq p-1\}.
$$
We have
$$
S_p=\{j\in\Z:\,-\floor{p/5}\leq j\leq\floor{p/5}\},
$$
and
$$
\widehat{S}_p=\cases
\{j\in\Z:\,-\floor{p/5}+1\leq j\leq\floor{p/5}\}&\text{if }p\equiv 1\mo{5},\\
\{j\in\Z:\,-\floor{p/5}\leq j\leq\floor{p/5}\}&\text{if }p\equiv 2,3\mo{5},\\
\{j\in\Z:\,-\floor{p/5}\leq j\leq\floor{p/5}+1\}&\text{if }p\equiv 4\mo{5}.\\
\endcases
$$
\endproclaim

These two lemmas above can be verified directly, so we omit the details here.

\noindent{\it Proof of Theorem 1.1}.
By (2.1), we have
$$
\cf_{p+1}(q)=\sum_{j=-\infty}^{\infty}(-1)^jq^{j(5j+1)/2}\qbinom{p}{\lfloor(p-5j)/2\rfloor}{q}.
$$
Observe that
$$
\qbinom{p}{k}{q}\equiv\cases1\pmod{[p]_q}\qquad&\text{if }k=1\text{ or }p,\\0\pmod{[p]_q}\qquad&\text{if }1\leq k\leq p-1.\endcases
$$
Then it follows that
$$
\cf_{p+1}(q)\equiv\sum_{\floor{(p-5j)/2}=0\text{ or }p}(-1)^jq^{j(5j+1)/2}\pmod{[p]_q}.
$$
It is easy to check that
$$
\{j:\,\floor{(p-5j)/2}=0\text{ or }p\}=\cases
\{(p-1)/5\}\qquad&\text{if }p\equiv1\mo{5},\\
\{-(p+1)/5\}\qquad&\text{if }p\equiv4\mo{5},\\
\emptyset\qquad&\text{if }p\equiv2,3\mo{5}.
\endcases
$$
Thus
$$
\cf_{p+1}(q)\equiv\sum_{\floor{(p-5j)/2}=0\text{ or }p}(-1)^jq^{j(5j+1)/2}\equiv\cases
1\pmod{[p]_q}\ \ &\text{if }p\equiv1,4\mo{5},\\
0\pmod{[p]_q}\ \ &\text{if }p\equiv2,3\mo{5}.
\endcases
$$
This concludes the proof of (1.4).

Also, applying (2.1) and Lemma 2.2, we deduce that
$$
\align
\cf_{p}(q)=&\sum_{0\leq\floor{(p-1-5j)/2}\leq p-1}(-1)^jq^{j(5j+1)/2}\qbinom{p-1}{\floor{(p-1-5j)/2}}{q}\\
=&\sum_{j=-\floor{p/5}}^{\floor{p/5}}(-1)^jq^{j(5j+1)/2}\prod_{k=1}^{\floor{(p-1-5j)/2}}\frac{[p-k]_q}{[k]_q}\\
=&\sum_{j=-\floor{p/5}}^{\floor{p/5}}(-1)^jq^{j(5j+1)/2}\prod_{k=1}^{\floor{(p-1-5j)/2}}\frac{[p]_q-[k]_q}{q^k[k]_q}\\
\equiv&\sum_{j=-\floor{p/5}}^{\floor{p/5}}(-1)^{j+\floor{(p-1-5j)/2}}q^{L(j)}\pmod{[p]_q}.
\endalign
$$
Suppose that $p\equiv 1\mo{5}$. Noting that $(p-1)/5$ is even and
$$
(2j-1)+\floor{(p-1-5(2j-1))/2}=2j+\floor{(p-1-5\cdot 2j)/2}+1=(p-1)/2-3j+1,
$$
we have
$$
\align
\cf_{p}(q)\equiv&(-1)^{-(p-1)/5}q^{L(-(p-1)/5)}+\sum_{j=-(p-6)/5}^{(p-1)/5}(-1)^{j+\lfloor(p-1-5j)/2\rfloor}q^{L(j)}\\
=&q^{(p-1)(p-2)/10-p(p-1)/2}+\sum_{j=-(p-11)/10}^{(p-1)/10}(-1)^{(p-1)/2-3j}(q^{L(2j)}-q^{L(2j-1)})\\
\equiv&q^{(4p+1)/5}\pmod{[p]_q}.
\endalign
$$
where Lemma 2.1 is applied in the last step. Similarly we obtain that
$$
\cf_{p}\equiv\cases
(-1)^{(p-2)/5}q^{L((p-2)/5)}\equiv-q^{(2p+1)/5}\pmod{[p]_q}\quad&\text{if }p\equiv 2\mo{5},\\
(-1)^{-(p-3)/5}q^{L(-(p-3)/5)}\equiv-q^{(3p+1)/5}\pmod{[p]_q}\quad&\text{if }p\equiv 3\mo{5},\\
(-1)^{(p-4)/5}q^{L((p-4)/5)}\equiv q^{(p+1)/5}\pmod{[p]_q}\quad&\text{if }p\equiv 4\mo{5}.\\
\endcases
$$
\qed

\noindent{\it Proof of Theorem 1.2}. In view of (2.2),
$$
\cg_{p}(q)=\sum_{j=-\infty}^{\infty}(-1)^jq^{j(5j-3)/2}\qbinom{p}{\lfloor(p-5j)/2\rfloor+1}{q}.
$$
And we have
$$
\{j:\,\lfloor(p-5j)/2\rfloor+1=0\text{ or }p\}=\cases
\{-(p-1)/5\}\qquad&\text{if }p\equiv1\mo{5},\\
\{-(p-2)/5\}\qquad&\text{if }p\equiv2\mo{5},\\
\{(p+2)/5\}\qquad&\text{if }p\equiv3\mo{5},\\
\{(p+1)/5\}\qquad&\text{if }p\equiv4\mo{5}.
\endcases
$$
Therefore
$$
\cg_{p}(q)\equiv\cases
(-1)^{-(p-1)/5}q^{(p+2)(p-1)/10}\equiv q^{(p-1)/5}\pmod{[p]_q}&\text{if }p\equiv1\mo{5},\\
(-1)^{-(p-2)/5}q^{(p+1)(p-2)/10}\equiv -q^{(3p-1)/5}\pmod{[p]_q}&\text{if }p\equiv2\mo{5},\\
(-1)^{(p+2)/5}q^{(p+2)(p-1)/10}\equiv -q^{(2p-1)/5}\pmod{[p]_q}&\text{if }p\equiv3\mo{5},\\
(-1)^{(p+1)/5}q^{(p+1)(p-2)/10}\equiv q^{(4p-1)/5}\pmod{[p]_q}&\text{if }p\equiv4\mo{5}.
\endcases
$$
Now let us turn to the proof of (1.7). From (2.2), it follows that
$$
\aligned
\cg_{p-1}(q)=&\sum_{j=-\infty}^{\infty}(-1)^jq^{j(5j-3)/2}\qbinom{p-1}{\lfloor(p-1-5j)/2\rfloor+1}{q}\\
=&\sum_{0\leq\lfloor(p-1-5j)/2\rfloor+1\leq p-1}^{\infty}(-1)^jq^{j(5j-3)/2}\qbinom{p-1}{\lfloor(p-1-5j)/2\rfloor+1}{q}.
\endaligned
$$
If $p\equiv 1\mo{5}$, then by Lemmas 2.1 and 2.2, we have
$$
\aligned
\cg_{p-1}(q)=&\sum_{j=-(p-1)/5+1}^{(p-1)/5}(-1)^jq^{j(5j-3)/2}\qbinom{p-1}{\lfloor(p-1-5j)/2\rfloor+1}{q}\\
\equiv&\sum_{j=-(p-6)/5}^{(p-1)/5}(-1)^{j+\lfloor(p-1-5j)/2\rfloor+1}q^{\widehat{L}(j)}\\
=&\sum_{j=-(p-11)/10}^{(p-1)/10}(-1)^{(p-1)/2-3j+1}(q^{\widehat{L}(2j)}-q^{\widehat{L}(2j-1)})\\
\equiv&0\pmod{[p]_q}.
\endaligned
$$
Similarly when $p\equiv 4\mo{5}$,
$$
\cg_{p-1}(q)\equiv\sum_{j=-(p-4)/5}^{(p+1)/5}(-1)^{j+\lfloor(p-1-5j)/2\rfloor+1}q^{\widehat{L}(j)}\equiv 0\pmod{[p]_q}.
$$
Finally, suppose that $p\equiv 2\text{ or }3\mo{5}$. Then
$$
\aligned
&\cg_{p-1}(q)\\
\equiv&\sum_{j=-\floor{p/5}}^{\floor{p/5}}(-1)^{j+\floor{(p-1-5j)/2}+1}q^{\widehat{L}(j)}\\
\equiv&\cases
(-1)^{(p-2)/5+1}q^{\widehat{L}((p-2)/5)}=q^{p(p-7)/10}\equiv 1\pmod{[p]_q}&\text{if }p\equiv 2\mo{5},\\
(-1)^{-(p-3)/5+(p-1)}q^{\widehat{L}(-(p-3)/5)}=q^{p(1-2p)/5}\equiv 1\pmod{[p]_q}&\text{if }p\equiv 3\mo{5}.
\endcases
\endaligned
$$
All are done.\qed
\heading
3. $q$-Pell number
\endheading
To prove Theorem 1.3, we need the similar identities as (2.1) and (2.2) for $\cp_n(q)$ and $\cq_n(q)$ respectively.
Fortunately, such identities have been established by Santos and Sills in [8]. Let
$$
\rT_1(n,m,q)=\sum_{j=0}^n(-q)^j\qbinom{n}{j}{q^2}\qbinom{2n-2j}{n-m-j}{q}.
$$
\proclaim{Lemma 3.1 {\rm (Santos and Sills, [8])}} Let $n$ be a non-negative integer. Then
$$
\align
\cp_{n+1}(q)=&\sum_{j=0}^n\sum_{k=0}^j q^{(j^2+j+k^2-k)/2}\qbinom{j}{k}{q}\qbinom{n-k}{j}{q}\\
=&\sum_{j=-\infty}^{\infty}(-1)^jq^{2j^2}\rT_1(n+1,4j+1,\sqrt{q}),\tag 3.1
\endalign
$$
and
$$
\align
\cq_{n+1}(q)=&\sum_{j=0}^n\sum_{k=0}^j q^{(j^2+j+k^2+k)/2}\qbinom{j}{k}{q}\qbinom{n-k}{j}{q}\\
=&\sum_{j=-\infty}^{\infty}(-1)^jq^{2j^2+j}\rT_1(n+1,4j+1,\sqrt{q}).\tag 3.2
\endalign
$$
\endproclaim

{\noindent{\it Proof of Theorem 1.3}}. For arbitrary two polynomials $f(q)$ and $g(q)$,
it is easy to see that $f(q^2)\mid g(q^2)$ implies $f(q)\mid g(q)$. Indeed, if
$$
g(q^2)/f(q^2)=\sum_{k=0}^n a_kq^k,
$$
then
$$
\sum_{k=0}^n a_k(-q)^k=g((-q)^2)/f((-q)^2)=g(q^2)/f(q^2)=\sum_{k=0}^n a_kq^k,
$$
whence $a_k=0$ for each odd $k$. Thus
$$
g(q)/f(q)=\sum_{0\leq k\leq n/2} a_{2k}q^k.
$$

Now it suffices to show that
$$
q^{(p^2-1)/4}\cp_p(q^2)\equiv\jacob{2}{p}\pmod{[p]_{q^2}}
$$
and
$$
\cq_p(q^2)\equiv\jacob{2}{p}q^{(p^2-1)/4}\pmod{[p]_{q^2}}.
$$
By (3.1), we have
$$
\align
\cp_{p}(q^2)=&\sum_{j=-\infty}^{\infty}(-1)^jq^{4j^2}\rT_1(n+1,4j+1,q)\\
=&\sum_{j=-\infty}^{\infty}(-1)^jq^{4j^2}\sum_{k=0}^{p}(-q)^{k}\qbinom{p}{k}{q^2}\qbinom{2p-2k}{p-4j-k-1}{q}\\
\equiv&\sum_{j=-\infty}^{\infty}(-1)^jq^{4j^2}\qbinom{2p}{p-4j-1}{q}\pmod{[p]_{q^2}}.
\endalign
$$
Since
$$
\qbinom{2p}{k}{q}\equiv 0\pmod{[p]_{q^2}}
$$
when $p\nmid k$, we have
$$
\align
\cp_{p}(q^2)\equiv&\cases
(-1)^{(p-1)/4}q^{(p-1)^2/4}\qquad&\text{if }p\equiv 1\mo{4}\\
(-1)^{-(p+1)/4}q^{(p+1)^2/4}\qquad&\text{if }p\equiv 3\mo{4}
\endcases\\
\equiv&\jacob{2}{p}q^{-(p^2-1)/4}\pmod{[p]_{q^2}}.
\endalign
$$
The proof of (1.10) is very similar. From (3.2), we deduce that
$$
\align
\cq_{p}(q^2)=&\sum_{j=-\infty}^{\infty}(-1)^jq^{4j^2+2j}\sum_{k=0}^{p}(-q)^{k}\qbinom{p}{k}{q^2}\qbinom{2p-2k}{p-4j-k-1}{q}\\
\equiv&\sum_{j=-\infty}^{\infty}(-1)^jq^{4j^2+2j}\qbinom{2p}{p-4j-1}{q}\\
\equiv&\cases
(-1)^{(p-1)/4}q^{(p-1)^2/4+(p-1)/2}\qquad&\text{if }p\equiv 1\mo{4}\\
(-1)^{-(p+1)/4}q^{(p+1)^2/4-(p+1)/2}\qquad&\text{if }p\equiv 3\mo{4}
\endcases\\
=&\jacob{2}{p}q^{(p^2-1)/4}\pmod{[p]_{q^2}}.
\endalign
$$
Thus the proofs of (1.9) and (1.10) are completed.

We need the following lemma in the proof of (1.11).
\proclaim{Lemma 3.2}
$$
\qbinom{2p+2}{p}{q}\equiv 1+q^p\pmod{[p]_{q^2}}
$$
for any odd prime $p$.
\endproclaim
\Proof. Write
$$
\qbinom{2p+2}{p}{q}=\frac{(1-q^{2p+2})(1-q^{2p+1})(1-q^{2p})\prod_{j=1}^{p-1}(1-q^{2p-j})}{(1-q^{p+2})(1-q^{p+1})(1-q^{p})\prod_{j=1}^{p-1}(1-q^{j})}.
$$
Since
$$
\frac{1-q^{2p-j}}{1-q^j}\equiv\frac{q^j-q^{2p}}{q^j(1-q^j)}=\frac{(1-q^{2p})-(1-q^j)}{q^j(1-q^j)}\equiv-q^{-j}\pmod{[p]_{q^2}},
$$
we have
$$
\qbinom{2p+2}{p}{q}\equiv(-1)^{p-1}q^{-\bi{p}{2}}\frac{(1-q^{2p+2})(1-q^{2p+1})(1-q^{2p})}{(1-q^{p+2})(1-q^{p+1})(1-q^{p})}\pmod{[p]_{q^2}}.
$$
And noting that
$$
q^{-\bi{p}{2}}\frac{(1-q^{2p+2})(1-q^{2p+1})}{(1-q^{p+2})(1-q^{p+1})}\equiv1\pmod{[p]_q},
$$
it follows that
$$
q^{-\bi{p}{2}}\frac{(1-q^{2p+2})(1-q^{2p+1})(1-q^{2p})}{(1-q^{p+2})(1-q^{p+1})(1-q^{p})}\equiv1+q^p\pmod{(1+q^p)[p]_q=(1+q)[p]_{q^2}}.
$$
We are done.\qed

From the definition of $\rT_1$, we have
$$
\align
&\rT_1(p+1,4j+1,q)+q\rT_1(p,4j+1,q)\\
=&\sum_{k=0}^{p+1}(-q)^k\qbinom{p+1}{k}{q^2}\qbinom{2p+2-2k}{p-4j-k}{q}-\sum_{k=0}^p(-q)^{k+1}\qbinom{p}{k}{q^2}\qbinom{2p-2k}{p-4j-k-1}{q}\\
=&\qbinom{2p+2}{p-4j}{q}+\sum_{k=0}^p(-q)^{k+1}\bigg(\qbinom{p+1}{k+1}{q^2}-\qbinom{p}{k}{q^2}\bigg)\qbinom{2p-2k}{p-4j-k-1}{q}\\
=&\qbinom{2p+2}{p-4j}{q}+\sum_{k=0}^p(-q)^{k+1}q^{{2(k+1)}}\qbinom{p}{k+1}{q^2}\qbinom{2p-2k}{p-4j-k-1}{q}\\
\equiv&\qbinom{2p+2}{p-4j}{q}+(-q)^p\qbinom{2}{-4j}{q}\pmod{[p]_{q^2}}.
\endalign
$$
Thus
$$
\align
\cp_{p+1}(q^2)+q\cp_{p}(q^2)
=&\sum_{j=-\infty}^{\infty}(-1)^jq^{4j^2}(\rT_1(p+1,4j+1,q)-q\rT_1(p,4j+1,q))\\
\equiv&\sum_{j=-\infty}^{\infty}(-1)^jq^{4j^2}\qbinom{2p+2}{p-4j}{q}
+\sum_{j=-\infty}^{\infty}(-1)^{p-j}q^{4j^2+p}\qbinom{2}{-4j}{q}\\
\equiv&\sum_{j=-\infty}^{\infty}(-1)^jq^{4j^2}\qbinom{2p+2}{p-4j}{q}-q^{p}\pmod{[p]_{q^2}}.
\endalign
$$
Note that
$$
\qbinom{2p+2}{k}{q}\not\equiv 0\pmod{[p]_{q^2}}
$$
only if $k\in\{0,1,2,p,p+1,p+2,2p,2p+1,2p+2\}$. Hence we have
$$
\align
&\sum_{j=-\infty}^{\infty}(-1)^jq^{4j^2}\qbinom{2p+2}{p-4j}{q}\\
\equiv&\qbinom{2p+2}{p}{q}+\cases
(-1)^{(p-1)/4}q^{(p-1)^2/4}\qbinom{2p+2}{1}{q}&\qquad\text{if }p\equiv 1\mo 4\\
(-1)^{-(p+1)/4}q^{(p+1)^2/4}\qbinom{2p+2}{2p+1}{q}&\qquad\text{if }p\equiv 3\mo 4
\endcases\\
\equiv&1+q^p+\jacob{2}{p}q^{-(p^2-1)/4}(1+q)\\
\equiv&1+q^p+(1+q)\cp_{p}(q^2)\pmod{[p]_{q^2}}.
\endalign
$$
This concludes that
$$
\cp_{p+1}(q^2)\equiv1+\cp_{p}(q^2)\pmod{[p]_{q^2}}.
$$
Similarly, we have
$$
\align
\cq_{p+1}(q^2)+q\cq_{p}(q^2)
=&\sum_{j=-\infty}^{\infty}(-1)^jq^{4j^2+2j}(\rT_1(p+1,4j+1,q)-q\rT_1(p,4j+1,q))\\
\equiv&\sum_{j=-\infty}^{\infty}(-1)^jq^{4j^2+2j}\qbinom{2p+2}{p-4j}{q}-q^{p}\\
\equiv&\qbinom{2p+2}{p}{q}+\jacob{2}{p}q^{(p^2-1)/4}(1+q)-q^{p}\\
\equiv&1+q^p+(1+q)\cq_{p}(q^2)-q^{p}\pmod{[p]_{q^2}},
\endalign
$$
whence
$$
\cq_{p+1}(q^2)\equiv1+\cq_{p}(q^2)\pmod{[p]_{q^2}}.
$$
\qed

\Ack. I thank my advisor, Zhi-Wei Sun, for his help on this paper.

\enddocument